\documentclass[12pt,oneside]{amsart}
\usepackage{bm, calc,geometry, verbatim, graphicx, amssymb, color, mathabx, mathtools, enumitem, bm, mathrsfs, amsmath, amsthm, ulem, xspace,tabularx,physics}
\usepackage[usenames,dvipsnames,svgnames,table]{xcolor}
\usepackage[pdftex,bookmarks,colorlinks,breaklinks]{hyperref} 

\hypersetup{linkcolor=blue,citecolor=red,filecolor=dullmagenta,urlcolor=darkblue} 
\newtheorem{theorem}{Theorem}[section]
\newtheorem{corollary}[theorem]{Corollary}
\newtheorem{lemma}[theorem]{Lemma}
\newtheorem{conjecture}[theorem]{Conjecture}

\newtheorem{problem}{Problem} 
\usepackage{physics}
\newcommand{\Mod}[1]{\ (\mathrm{mod}\ #1)}

\newcommand\mult{\operatorname{\textup{{\fontfamily{ptm}\selectfont mult}}}}
\newcommand\dg{\operatorname{\textup{{\fontfamily{ptm}\selectfont deg}}}}

\newcommand{\cc}{\mathcal}

\newcommand\roundup[1]{\left\lceil#1\right\rceil}
\newcommand\rounddown[1]{\left\lfloor#1\right\rfloor}

     \makeatletter
      \def\@setcopyright{}
      \def\serieslogo@{}
      \makeatother      
\begin{document}
   \author{Amin  Bahmanian}
   \address{Department of Mathematics,
  Illinois State University, Normal, IL USA 61790-4520}
\title[Toward a Three-dimensional Counterpart of Cruse's Theorem]{Toward a Three-dimensional Counterpart of Cruse's Theorem}

   \begin{abstract}  
Completing partial latin squares is NP-complete. Motivated by Ryser's theorem for latin rectangles, in 1974, Cruse found conditions that ensure a partial symmetric latin square of order $m$ can be embedded in a symmetric latin square of order $n$. Loosely speaking, this results asserts that an $n$-coloring of the edges of the complete $m$-vertex graph $K_m$  can be embedded in a one-factorization of $K_n$ if and only if $n$ is even and the number of edges of each color is at least $m-n/2$.
We establish necessary and sufficient conditions under which an edge-coloring of the complete $\lambda$-fold $m$-vertex 3-graph $\lambda K_m^3$  can be embedded  in a one-factorization of $\lambda K_n^3$. In particular, we prove the first known Ryser type theorem for hypergraphs by showing that if $n\equiv 0\Mod 3$,  any  edge-coloring of  $\lambda K_m^3$  where the number of triples of each color is at least $m/2-n/6$,  can be embedded in a one-factorization of $\lambda K_n^3$. Finally we prove an Evans type result by showing that if $n\equiv 0\Mod 3$ and $n\geq 3m$, then any $q$-coloring of the edges of any $F\subseteq\lambda K_m^3$ can be embedded in a one-factorization of $\lambda K_n^3$  as long as $q\leq \lambda \binom{n-1}{2}-\lambda \binom{m}{3}/\rounddown{m/3}$. 

   \end{abstract}
   \date{\today}
   \maketitle   
\section{Introduction}   

Ryser showed that any $r\times s$ latin rectangle $L$ can be embedded in an $n\times n$ latin square if and only if the number of occurrences of each symbol in $L$ is at least $r+s-n$ \cite{MR42361}. In other words, an $n$-coloring of the edges of the complete bipartite graph $K_{r,s}$ can be embedded in an $n$-coloring of $K_{n,n}$ if and only if the number of edges of each color in  $K_{r,s}$ is at least $r+s-n$. The main motivation of this note is to find a three-dimensional analogue of Ryser's theorem for symmetric latin cubes  for which there has been very little success over the last seventy years. To our best knowledge, our main result is the first result of this kind. The study of latin cubes was initiated in the 1940's by Fisher \cite{MR13113} and   Kishen \cite {Kishen42,MR34743}. These objects  are also closely related to orthogonal arrays \cite{MR1693498, MR3495977} which themselves  have  applications in  statistics, coding theory, and cryptography \cite{MR2246267}.  

Let $L$ be an $n\times n\times n$ array. We imagine $L$ as a 3-dimensional array having layers stacked on top of each other, and we call it a cube. The symbol in position $(i,j,k)$ of $L$ is denoted by $L_{ijk}$. A {\it layer} in $L$ is obtained by fixing one coordinate. Hence, each layer is parallel to one of the faces of the cube $L$.  We say that $L$ is a {\it layer-equitable latin cube} if it is filled with $n$ different symbols such that each layer contains each symbol $n$ times. Stacking 
\begin{tabular}{ c c c }
 1 & 1 & 2  \\ 
 1 & 2 & 3 \\
 2 & 3 & 3
\end{tabular}
on top of  
\begin{tabular}{ c c c }
 2 & 2 & 3  \\ 
 2 & 3 & 1 \\
 3 & 1 & 1
\end{tabular}
on top of
\begin{tabular}{ c c c }
 3 & 3 & 1  \\ 
 3 & 1 & 2 \\
 1 & 2 & 2
\end{tabular}
forms a layer-equitable latin cube of order 3. A layer-equitable latin cube in which each layer is a latin square, is  called a {\it layer-latin latin cube}. Stacking 
\begin{tabular}{ c c }
 1 & 2  \\ 
 2 & 1     
\end{tabular}
on top of  
\begin{tabular}{ c c }
 2 & 1  \\ 
 1 & 2     
\end{tabular}
forms a layer-latin latin cube of order 2. If $L$ is filled with $n^2$ symbols such that in each layer every symbol occurs exactly once, then $L$ is a {\it layer-rainbow latin cube}. Stacking 
\begin{tabular}{ c c }
 1 & 2  \\ 
 3 & 4     
\end{tabular}
on top of  
\begin{tabular}{ c c }
 4 & 3  \\ 
 2 & 1     
\end{tabular}
forms a layer-rainbow latin cube of order 2. There is no consensus on which of the above three notions should be called a latin cube, but most statisticians refer to layer-rainbow latin cubes as latin cubes, and most combinatorialists refer to  layer-latin latin cubes as permutation cubes or latin cubes. For the sake of convenience, throughout the rest of this note we shall refer to layer-latin latin cubes as latin cubes.

Ryser's theorem implies that every $r\times n$ latin rectangle on $n$ symbols can be embedded in an $n\times n$ latin square \cite{MR13111}. This result is not true in the 3-dimensional case, as there are many examples of $r\times n\times n$ partial latin cubes that cannot be completed to $n\times n\times n$ latin cubes \cite{FuNotinMath, MR1352785, MR2399374}. There are only a few results on embedding  partial latin cubes, and in fact, most of the current literature is focused on the cases where a partial latin cube cannot be completed. 
Cruse showed that  a partial  latin cube of order $n$ can be embedded in a  latin cube of order $16n^4$ \cite{MR345854}; this result was extended to idempotent  latin cubes by Lindner \cite{MR441759} and Csima \cite{MR539421}. Potapov  improved Cruse's bound to  $n^3$ \cite{MR4453977}. 
Denley and \"{O}hman \cite{MR3186483} found sufficient conditions for when certain $k \times l \times  m$ partial Latin cubes, namely latin boxes (or latin parallelepipeds), can be embedded in $k \times n \times m$ latin boxes, and subsequently, to $k \times n \times n$ latin boxes, though they were  not  able to extend in the third dimension to obtain  latin cubes. For related embedding results in higher dimensions and connections with maximum distance separable  codes, see Krotov and Sotnikova \cite{MR3357771}, and Potapov \cite{MR2961772, MR4453977}. Various results on latin boxes that cannot be completed to latin cubes can be found in \cite{ MR3336595, MR2902643, MR872667, MR677575, MR1094977, MR1018253, MR1352785, MR2979500, MR2838020}. For further results on latin cubes, see \cite{MR4469221, MR3897549, MR3600882, MR2246267LaywMull, MR4113541} and references therein, and for recent results on the number of factorizations, see \cite{2017arXiv170505225L}
 and \cite{MR3603556}. 

 Cruse settled conditions that ensure an $m\times m$ symmetric latin rectangle on $n$ symbols can be embedded in a symmetric latin square of order $n$ \cite{MR329925} (see Theorem \ref{crusorigthmgt}). We consider the following three-dimensional analogue of Cruse's theorem. 
 \begin{problem} \label{mainprob}
     Find conditions under which an $m\times m\times m$ symmetric layer-rainbow latin box on $n$ symbols can be embedded in a symmetric layer-rainbow latin cube of order $n$. 
 \end{problem}
To investigate this problem, first we need to define symmetry. A layer-rainbow latin cube $L$ of order $n$ is {\it symmetric} if  it meets the following conditions.
\begin{enumerate}
    \item $L_{ij\ell}=L_{j\ell i}=L_{\ell ij}$ for distinct $i,j,\ell \in  \{1,\dots,n\}$,
    \item  $L_{iij}=L_{jj i}, L_{iji}=L_{jij}, L_{ij j}=L_{jii}$ for  $i,j\in \{1,\dots,n\}$.
\end{enumerate}  
Here is an example of a   symmetric layer-rainbow latin cube of order 5. 
\begin{center} 
\begin{tabular}{ |c |c| c| c| c|}\hline
A & F & K & P & U \\ \hline
B & G & L & Q & V \\ \hline
C & H & M & R & W \\ \hline
D & I & N & S & X \\ \hline
E & J & O & T & Y \\ \hline
\end{tabular}\quad 
\begin{tabular}{ |c |c| c| c| c|}\hline
G & B & H & I & J \\ \hline
F & M & X & W & P \\ \hline
L & T & U & Y & S \\ \hline
Q & O & E & A & K \\ \hline
V & R & D & C & N \\ \hline
\end{tabular}\quad 
\begin{tabular}{ |c |c| c| c| c|}\hline
M & L & C & N & O \\ \hline
H & U & T & E & D \\ \hline
K & X & P & J & A \\ \hline
R & Y & V & F & B \\ \hline
W & S & Q & G & I \\ \hline
\end{tabular}
\end{center}
\begin{center}
\begin{tabular}{ |c |c| c| c| c|}\hline
S & Q & R & D & T \\ \hline
I & A & Y & O & C \\ \hline
N & E & F & V & G \\ \hline
P & W & J & U & H \\ \hline
X & K & B & L & M \\ \hline
\end{tabular}\quad 
\begin{tabular}{ |c |c| c| c| c|}\hline
Y & V & W & X & E \\ \hline
J & N & S & K & R \\ \hline
O & D & I & B & Q \\ \hline
T & C & G & M & L \\ \hline
U & P & A & H & F \\ \hline
\end{tabular} 
\end{center}
For convenience, we have used letters for symbols. For example, $L_{112}=L_{221}=F$, $L_{255}=L_{522}=N$, and $L_{134}=L_{341}=L_{413}=R$. These objects are natural three-dimensional analogues of symmetric latin squares. Note that in a symmetric latin square, if row $i$ is given, column $i$ is uniquely determined (if fact row $i$ and column $i$ are identical). For $i\in \{1,\dots,n\}$, let 
\begin{align*}
    L_{i**}=\{(i,j,\ell, L_{ij\ell})\ |\  1\leq j,\ell\leq n\},\\
    L_{*i*}=\{(j,i,\ell, L_{ji\ell})\ |\  1\leq j,\ell\leq n\},\\
    L_{**i}=\{(j,\ell, i,L_{j\ell i})\ |\  1\leq j,\ell\leq n\}.
\end{align*}
In a symmetric layer-rainbow latin cube, if one of the three layers in $\{L_{i**},L_{*i*},L_{**i} \}$ is given, the other two layers can be uniquely determined. To see this, suppose $L_{i**}$ is given. To  determine $L_{*i*}$, observe that for distinct $j, \ell \in [n]$, we have  $L_{iij}\subseteq L_{i**}$, $L_{jii}=L_{ijj}$ and $L_{ji\ell}=L_{i\ell j}$,  . The argument for other cases is  similar.

For $k\in \mathbb N$, and a finite set $X$ of {\it vertices}, let  $[k]=\{1,\dots,k\}$, and let $K_{ m}^k$ be the collection of all subsets of size $k$ of $X$.  Let $E\subseteq \mathcal P(X)$, where $\mathcal P(X)$ is the power set of $X$. We call $G=(X,E)$ a {\it simple hypergraph} and each element of $E$ is an {\it edge}. Whenever it is not ambiguous, we identify a hypergraph with its edge set. In a {\it general hypergraph}, we allow  edges to be repeated and in addition we allow the vertices to be repeated within each edge. We use $\mult_G(e)$ to denote the number of occurrences of $e$ in $G$. The {\it degree} of  $x\in X$ in $G$, written $\dg_G(x)$, is the number of occurrences of $x$ in $G$. For $\lambda\in \mathbb N$,  $\lambda G$ is obtained by replacing each $e\in E$ by $\lambda$ copies of $e$. We call $\lambda K_n^k$ {\it the complete $\lambda$-fold $n$-vertex $k$-graph}, and we abbreviate $K_n^2$ to $K_n$. 
A {\it $k$-coloring} of $G$ is  a partition of $E$ into {\it color classes} $G(1),\dots,G(k)$ such that in each color class each pair of edges are disjoint (so we are only concerned with edge-colorings). If each color class is a partition of $X$ (i.e. a perfect matching), then the coloring is a {\it one-factorization}.

We restate    Cruse's theorem   as follows.
\begin{theorem} \cite[Theorem 1]{MR329925} \label{crusorigthmgt}
An $n$-coloring of  $F:=K_m$ can be embedded in a one-factorization of $K_n$  if and only if $n$ is even and
$$|F(i)|\geq m-\frac{n}{2}   \quad \mbox { for } i\in [n].$$
\end{theorem}
Symmetric latin squares of order $n$ can be seen as one-factorizations of $K_n^1\cup K_n^2$, so to account for  $K_n^1$, a further condition is needed (see \cite{MR329925}). Symmetric layer-rainbow latin cubes  can be seen as one-factorizations of $K_n^1\cup 3K_n^2\cup 2K_n^3$ \cite{doi:10.1137/22M1494488}.
For example, consider the following one-factorization $\cc F$ of $K_5^1\cup 3K_5^2\cup 2K_5^3$ (We abbreviate sets such as $\{x,y,z\}$ to $xyz$). 
\begin{align*}
\big\{&\{1,24, 35\},\{ 12, 345\}, \{13, 245\}, \{ 14, 235\},\{ 15, 234\}, \{ 5,12, 34\}, \{12, 345\}, \{45,123\}, \\
&\{35, 124\}, \{34, 125\}, \{13, 245\}, \{45, 123\}, \{2, 13,45\}, \{25, 134\}, \{24, 135\}, \{3,14, 25\}, \\
&\{35, 124\}, \{25, 134\}, \{14, 235\}, \{23, 145\}, \{4, 15, 23\}, \{34, 125\}, \{24, 135\}, \\
&\{23, 145\}, \{15, 234\}\big\}.
\end{align*}
Let us assign symbols $A, B,\dots, Y$ to the 25 color classes, respectively. Then we fill $L_{iii}$ with the color of  $i$ in $\cc F$ for  $i\in \{1,\dots, 5\}$. For distinct  $i,j\in \{1,\dots, 5\}$, suppose the color of the three copies of $\{ij\}$ in $\cc F$ are $c,c'$ and $c''$. We fill  the cells $L_{iij}$ and $L_{jji}$ with $c$, $L_{iji}$ and $L_{jij}$ with $c'$, and $L_{jii}$ and $L_{ijj}$ with $c''$. For distinct  $i,j, \ell\in \{1,\dots, 5\}$, suppose the color of the two copies of $\{ij \ell\}$ are $c$ and $c'$. We fill the cells  $L_{ij\ell}$, $L_{j\ell i}$ and $L_{\ell ij}$ with $c$, and $L_{i\ell j}$, $L_{\ell j i}$, and $L_{ji \ell}$ with $c'$. The result is the symmetric layer-rainbow latin cube of order 5 we mentioned earlier. We can obtain other similar symmetric layer-rainbow latin  cubes using $\cc F$.

Now, we are ready to restate Problem \ref{mainprob} more precisely.
\begin{problem} \label{mainprobprecise}
    Find conditions under which an $n^2$-coloring of $K_m^1\cup 3K_m^2\cup 2K_m^3$ can be embedded in a one-factorization of $K_n^1\cup 3K_n^2\cup 2K_n^3$.
\end{problem}

Here, we establish conditions that ensure a coloring of $\lambda K_m^3$ can be embedded in a one-factorization of $\lambda K_n^3$ (see Theorem \ref{embedvialistthm}). 
\begin{theorem} \label{weakerthm}
Let $k=\lambda \binom{n-1}{2}$. A $k$-coloring of $F:=\lambda K_m^3$ satisfying 
\begin{align} \label{weakercond}
    |F(i)|\geq \frac{m}{2}-\frac{n}{6}\quad \mbox { for } i\in [k],
\end{align}
can be embedded in a one-factorization of $\lambda K_n^3$ if and only if $n\equiv 0\Mod 3$.
\end{theorem}
We remark that for $m\equiv n\equiv 0 \Mod 3$ and $n\leq 2m-1$, there exist  $\binom{n-1}{2}$-colorings of $K_m^3$ which cannot  be embedded in a one-factorization of $K_n^3$ \cite{MR1249714}.
Evans, and independently, S. K. Stein \cite{MR122728} used Ryser's theorem to show that a partial latin square of  order $m$ on $m$ symbols can be embedded in a latin square of order $n$ for $n\geq 2m$. Cruse extended this idea to symmetric latin squares which we restate as follows. 
\begin{theorem} \cite{MR329925}
If  $ n\geq 2m$, then an $(m-1)$-coloring of any $F \subseteq  K_m$ can be embedded in a one-factorization  of $K_n$ if and only if $ n\equiv 0 \Mod 2$.
\end{theorem}
Here is our next result. 
\begin{theorem} \label{evanstype} 
If  $ n\geq 3m$ and 
$$q\leq \lambda \binom{ n-1}{2}-\dfrac{\lambda \dbinom{m}{3}}{\rounddown{\dfrac{m}{3}}},$$ 
then a $q$-coloring of any $F \subseteq \lambda  K_m^3$ can be embedded in a one-factorization  of $\lambda  K_n^3$ if and only if $ n\equiv 0 \Mod 3$.
\end{theorem}
Note that for $n\geq 3m$, $\binom{ n-1}{2}- \binom{m}{3}/\rounddown{m/3}\geq \binom{m-1}{2}$. 

Theorems \ref{weakerthm}, \ref{evanstype}, and \ref{embedvialistthm}  can be restated as results on embedding partial  symmetric layer-rainbow latin cubes in  partial symmetric layer-rainbow latin cubes where all diagonal entries are empty (Here, diagonal entries of a cube $L$ are those $L_{ijk}$s where $|\{i,j,k\}|\neq 3$).

The proof of Theorem \ref{weakerthm} relies on a stronger result, Theorem \ref{embedvialistthm}, together with a list-chromatic index bound for  complete graphs due to H\"{a}ggkvist and Janssen \cite{MR1464567}.  For a function $\gamma$ from the edges of $G$ to $\mathcal P([k])$, a $\gamma$-coloring of $G$ is a coloring of $G$ in which the color of each edge $e$ of  $G$ is chosen from the set $\gamma(e)$ and that each color class is a matching. Here is our main result. 
\begin{theorem} \label{embedvialistthm}
Let $F= \lambda  K_m^3,  G= \lambda  K_n^3, H=\lambda ( n-m)K_m,k=\lambda \binom{ n-1}{2}$, and let $X$ be the vertex set of $F$. 
A $k$-coloring of $F$ can be embedded in a one-factorization of $G$ if and only if $ n\equiv 0 \Mod 3$, and  there exists a  $\gamma$-coloring of $H$ such that 
\begin{align} \label{mainlistcond}
    |H(i)|\geq m-\frac{ n}{3}-2|F(i)|\quad \mbox { for } i\in [k],
\end{align}
where
\begin{align} \label{gammacond}
\gamma(uv)=\{ i\in [k]\ |\ \dg_{F(i)}(u)=\dg_{F(i)}(v)=0\} \quad \mbox { for } u,v\in X.
\end{align}
\end{theorem}
We remark that  Baranyai \cite{MR0416986} showed that $K_n^k$ is one-factorable if and only if $n\equiv 0\Mod k$. 

The paper is organized as follows. In Section \ref{prereq}, we provide the prerequisites. We prove our main result in Section \ref{mainressec}, followed by the proof of Theorem \ref{weakerthm} in Section \ref{weakerthmsec}. Finally, in Section \ref{remarksec}, we consider our  main problem without assuming that $ n\equiv 0\Mod 3$, and propose a few problems. To avoid trivial cases, we will always assume that $m\geq 3$ and $ n>m$. We abbreviate sets like $\{u,v\}$ to $uv$,  $\{u,v,w\}$ to $uvw$,  $\{u,u,v\}$ to $u^2v$,   $\{u,u,u\}$ to $u^3$, etc..

\section{List Colorings and Detachments} \label{prereq}
Recall that for a function $\gamma$ from the edges of $G$ to $\mathcal P([k])$, a $\gamma$-coloring of $G$ is a coloring of $G$ in which the color of each edge $e$ of  $G$ is chosen from the set $\gamma(e)$. H\"{a}ggkvist and Janssen \cite{MR1464567} showed that $K_m$ has a $\gamma$-coloring if  $|\gamma(e)|\geq m$ for each edge $e$ of $K_m$.  Consequently, 
\begin{align} \label{haggjenmult}
   \lambda K_m \mbox { has a } \gamma\mbox{-coloring}\quad  \mbox { if } |\gamma(e)|\geq \lambda m \mbox{ for each edge }e \mbox{ of } \lambda K_m.
\end{align}

Throughout this paper, $a\approx b$ means $a\in \{\rounddown{b}, \roundup{b}\}$. In a hypergraph an edge may  occur multiple times, and  a vertex may occur  multiple times within an edge. In other words, the edge set and every element of the edge set is a multiset. By an edge of the form $\{u_1^{m_1},\dots,u_s^{m_s}\}$, we mean an edge in which vertex $u_i$ occurs $m_i$ times for $1\leq i\leq s$.

Let $p$ be a positive integer. Given a hypergraph $G$ with vertex set $X$ whose edge set $E$ is colored, and $\alpha\in X$, we obtain a hypergraph $F$ on $|X|+p-1$ vertices and $|E|$ edges in the following way.
\begin{enumerate} 
\item [(a)] Replace $\alpha$ in $G$ by $p$ new vertices $\alpha_1,\dots,\alpha_p$ in $F$;
\item [(b)] Replace each edge $\{\alpha^q\}\cup U$ in $G$ by an edge $\{\alpha_1^{q_1},\dots,\alpha_p^{q_p}\}\cup U$ in $F$, where $\alpha\notin U\subseteq X$ and $q_1+\dots+q_p=q\geq 1$. 
\item [(c)] Leave the remaining vertices and edges of $G$ intact and add them to $F$.
\end{enumerate} 
The resulting hypergraph $F$ is called an {\it $(\alpha,p)$-detachment}  of $G$, and $G$ is an {\it amalgamation} of $F$ obtained by {\it identifying} $\alpha_1,\dots, \alpha_p$ by $\alpha$.  It is clear that $F$ is in general far from being uniquely determined. We need the following detachment lemma. 
\begin{lemma}  \cite[Theorem 4.1]{MR2942724}
\label{amalgambahrod} 
  Let $G$ be a hypergraph with vertex set $X$ whose edges are colored with $k$  colors, and let $\alpha$ be a vertex of $G$. There exists an $(\alpha,p)$-detachment $F$ of $G$ obtained by splitting $\alpha$ in $G$ into $\alpha_1,\dots,\alpha_p$ in $F$, such that the following conditions are satisfied.
\begin{enumerate}
     \item [\textup{(i)}]  $\dg_{F(j)}(\alpha_i)\approx \dg_{G(j)}(\alpha)/p$  for  $i\in [p],j\in [k]$;
     \item [\textup{(ii)}]  $\mult_F(\alpha_i u v)\approx \mult_G(\alpha uv)/p$ for  $i\in[p]$ and distinct $u,v\in X\backslash \{\alpha\}$;
     \item [\textup{(iii)}]  $\mult_F(\alpha_i \alpha_j u)\approx \mult_G(\alpha^2 u)/\binom{p}{2}$ for distinct $i,j\in[p],u\in X\backslash \{\alpha\}$;
     \item [\textup{(iv)}]  $\mult_F(\alpha_i \alpha_j \alpha_{\ell})\approx \mult_G(\alpha^3 )/\binom{p}{3}$ for  distinct $i,j,\ell \in[p]$.
\end{enumerate}  
\end{lemma}

\section{Proof of Theorem \ref{embedvialistthm}} \label{mainressec}
To prove the necessity, suppose that a $k$-coloring of $F$ is embedded in a one-factorization  of $G$. The existence of a one-factor in $G$ implies that $ n\equiv 0 \Mod 3$. Let 
$$
H'=\{ e\cap X\ \big|\ e\in  G, |e\cap X|= 2\}.
$$
The $k$-coloring of $G$  induces a $\gamma$-coloring of $H'\cong H$. For $i\in [k]$, let $c_i$ and $d_i$ be the number of edges of $G$ colored $i$ with exactly 1 and 0 vertices, respectively, in $X$. We have
\begin{align*}
\begin{cases}
|F(i)|+|H(i)|+c_i+d_i=\dfrac{ n}{3} & \mbox { for } i\in [k],\\
   3|F(i)|+2|H(i)|+c_i=  m & \mbox { for } i\in [k].
\end{cases}
\end{align*} 
Therefore,
$$d_i= 2|F(i)|+|H(i)|-m+\dfrac{ n}{3}\quad \mbox { for } i\in [k].$$
Since $d_i\geq 0$ for $i\in [k]$, we have 
\begin{align*} 
    |H(i)|\geq m-\dfrac{ n}{3}-2|F(i)|\quad \mbox { for } i\in [k].
\end{align*}

To prove the sufficiency, suppose that a $k$-coloring of $F$ is given and that there exists a  $\gamma$-coloring of $H$ such that \eqref{mainlistcond} holds.
Let $G_1$ be obtained by adding a new vertex $\alpha$ and the following multiset 

$$\Big\{\alpha u v \ |\ u,v\in X, u\neq v\Big\}^{\lambda ( n-m)}$$ 
of edges to $F$. In other words, for each pair of distinct vertices $u,v\in X$, we add $\lambda ( n-m)$ copies of an edge of the form $\alpha u v$ to $G_1$. 
We use the $\gamma$-coloring of each edge $uv$ of $H$ to color the corresponding edge $\alpha u v$ of $G_1$. Observe that for $i\in [k]$, $$\sum_{u\in X}  \dg_{G_1(i)}(u)=3|F(i)|+2|H(i)|.$$

Let $G_2$ be obtained by adding the following multiset 
$$\Big\{\alpha^2 u \ |\ u\in X\Big\}^{\lambda \binom{ n-m}{2}}\bigcup \Big \{\alpha^3\Big\}^{\lambda\binom{ n-m}{3}}$$
of edges to $G_1$. In other words, we add $\lambda \binom{ n-m}{3}$ copies of an edge of the form $\alpha^3$, and for each  $u\in X$, we add $\lambda \binom{ n-m}{2}$ copies of an edge of the form $\alpha^2 u$ to $G_2$. Recall that by \eqref{mainlistcond}, $2|F(i)|+|H(i)|-m+\dfrac{ n}{3}$ is a non-negative integer for $i\in [k]$. We color $\alpha^2 u$-edges   (for $u\in X$), and then we color the $\alpha^3$-edges of $G_2$ such that 
\begin{align*}
\begin{cases}
\mult_{G_2(i)}(\alpha^2 u)= 1-\dg_{G_1(i)}(u)&\mbox { for } i\in [k],\\
\mult_{G_2(i)}(\alpha^3)= 2|F(i)|+|H(i)|-m+\dfrac{ n}{3}&\mbox { for } i\in [k].
\end{cases}
\end{align*}
This is possible, for
\begin{align*}
    \sum_{i\in [k]}\Big(1-\dg_{G_1(i)}(u)\Big)&=k-\dg_{G_1}(u)\\
    &= \lambda \binom{ n-1}{2}-\lambda \binom{m-1}{2}-\lambda (m-1)( n-m)\\
    &=\lambda \binom{ n-m}{2}\quad\quad\quad\quad\quad \quad\quad\quad\quad\quad \mbox { for } u\in X,
\end{align*}
and \begin{align*}
    \frac{1}{\lambda }\sum_{i\in [k]}& \Big(2|F(i)|+|H(i)|-m+\frac{ n}{3}\Big)\\
    &=2\binom{m}{3}+\left( n-m\right)\binom{m}{2}+\left(\frac{ n}{3}-m\right)\binom{ n-1}{2}\\
    &=\binom{ n-m}{3}.
\end{align*}

For $i\in [k]$, we have the following.
\begin{align*}
  \dg_{G_2(i)}(\alpha)&=3\mult_{G_2(i)}(\alpha^3)+2\sum_{u\in X} \mult_{G_2(i)}(\alpha^2 u)+|H(i)|\\
  &=\big(6|F(i)|+3|H(i)|-3m+ n\big)+2\sum_{u\in X}  \Big(1-\dg_{G_1(i)}(u)\Big)+|H(i)|\\
  &=6|F(i)|+4|H(i)|-3m+ n+\big(2m-6|F(i)|-4|H(i)|\big)\\
  &=  n-m.
\end{align*}
By the detachment lemma, there exists a hypergraph $G$,  obtained by replacing the vertex $\alpha$ of $G_2$ by $ n-m$ new vertices  $\alpha_1,\dots,\alpha_{ n-m}$ in $G$, replacing each $\alpha uv$-edge in $G_2$ by an $\alpha_i uv$-edge in $G$ (for some $i\in[ n-m]$), replacing each $\alpha^2 u$-edge in $G_2$ by an $\alpha_i \alpha_j u$-edge in $G$ (for distinct $i,j\in[ n-m]$), replacing each $\alpha^3$-edge in $G_2$ by an $\alpha_i \alpha_j \alpha_\ell$-edge in $G$ (for distinct $i,j, \ell\in[ n-m]$), such that the edges incident with in $G_2$ are shared as evenly as possible among $\alpha_1,\dots,\alpha_{ n-m}$ in $G$ in the following way.
\begin{itemize}
    \item [(a)] For  $i\in [ n-m]$, and $j\in[k]$
    \begin{align*} 
        \dg_{G(j)}(\alpha_i)=\dfrac{\dg_{G_2(j)}(\alpha)}{ n-m}=1;
    \end{align*}
    \item [(b)] For $i\in[ n-m]$ and distinct $u,v\in X$, 
    \begin{align*} 
    \mult_{G}(\alpha_i uv)=\dfrac{\mult_{G_2}(\alpha uv)}{ n-m}=\lambda;
       \end{align*}
    \item [(c)] For distinct $i_1,i_2\in[ n-m]$ and $u\in X$, 
    \begin{align*} 
    \mult_{G}(\alpha_{i_1} \alpha_{i_2} u)=\dfrac{\mult_{G_2}(\alpha^2 u)}{\binom{ n-m}{2}}=\lambda;
       \end{align*}
    \item [(d)] For distinct $i_1,i_2,i_3\in[ n-m]$, 
    \begin{align*} 
    \mult_{G}(\alpha_{i_1} \alpha_{i_2}\alpha_{i_3})=\dfrac{\mult_{G_2}(\alpha ^3)}{\binom{ n-m}{3}}=\lambda.
   \end{align*}
\end{itemize}
By (b)--(d), $G \cong \lambda  K_n^3$, and by (a), each color class of $G$ is a perfect matching. This completes the proof. \qed

\section{Proof of Theorem \ref{weakerthm}} \label{weakerthmsec}
Since \eqref{weakercond} holds,  \eqref{mainlistcond} is trivial, and so by Theorem \ref{embedvialistthm}, it is enough to show that $H:=\lambda ( n-m)K_m$ has a $\gamma$-coloring where $\gamma$ satisfies \eqref{gammacond}. Thus, by \eqref{haggjenmult}, if we show that $|\gamma(e)| \geq \lambda m( n-m)$ for each $e$ in  $H$, we are done. Let 
$$
C(u)=\{i\in [k] \ |\ \dg_{F(i)}(u)=1\} \quad \mbox { for } u\in X.
$$
Observe that  $|C(u)|=\dg_F(u)=\lambda \binom{m-1}{2}$ for  $u\in X$. Moreover, for distinct  $u,v\in X$,  there are exactly $\lambda (m-2)$ edges in $F$ containing both $u$ and $v$, and so $|C(u)\cap C(v)|\geq \lambda (m-2)$. 
For $e=uv\in H$, we have the following which completes the proof. 
\begin{align*}
    |\gamma(uv)| - \lambda m( n-m) &=|\{ i\in [k]\ |\ \dg_{F(i)}(u)=\dg_{F(i)}(v)=0\}|- \lambda m( n-m)\\
    &=|\{i\in [k] \ |\ i\in \overline{C(u)}\cap \overline{C(v)}\}|- \lambda m( n- m)\\
    &=|\overline{C(u)\cup C(v)}|- \lambda m( n- m)\\
    &=k-|C(u)\cup C(v)|-\lambda m( n-m)\\
    &=k-|C(u)|-|C(v)|+ |C(u)\cap C(v)|-\lambda m( n- m)\\
    &\geq \lambda\binom{ n-1}{2}-2\lambda\binom{ m-1}{2}+\lambda( m-2)-\lambda m( n- m)\\
    &=\frac{\lambda}{2}\left ( n^2 -2 m n- 3 n + 8  m-6   \right)\\
    &\geq \frac{\lambda}{2}\left( n^2-2 m n-3 n-2 m^2+6 m+2\right)\\
    &=\frac{6\lambda}{ n- m}  \binom{ m}{3}-\frac{\lambda}{ n- m}  \binom{ n-1}{2}\left(3  m- n\right)\\
    &=\frac{6}{ n- m}\sum_{i\in [k]} |F(i)|- \frac{6\lambda}{ n- m} \binom{ n-1}{2}\left(\frac{ m}{2}-\frac{ n}{6}\right) \geq 0.
\end{align*} \qed

An immediate consequence of Theorem \ref{weakerthm} is the following which improves the best previous bound of $ n\geq (2+\sqrt{2}) m$ \cite{MR3056885}. 
\begin{corollary} \label{n3mcor}
If $ n\geq 3 m$, then a $\lambda \binom{ n-1}{2}$-coloring of $\lambda  K_m^3$ can be embedded in a one-factorization  of $\lambda  K_n^3$ if and only if $ n\equiv 0\Mod 3$.
\end{corollary}

\section{Proof of Theorem \ref{evanstype}}
Let $k=\lambda\binom{ n-1}{2}$. Suppose that  $ n\geq 3 m$, $ n\equiv 0 \Mod 3$, and $k-q\geq \lambda \binom{ m}{3}/\rounddown{ m/3}$. We show that we can embed the given  $q$-coloring of  $F \subseteq \lambda  K_m^3$ in a $k$-coloring  of $\lambda  K_n^3$. Let $G$ be the hypergraph whose edge set is $\lambda \binom{ m}{3}$ copies of $\{\alpha^3\}$ and color $G$  with colors $\{q+1,\dots, k\}$ such that 
$$\mult_{G(i)}(\alpha^3)\leq\rounddown{\frac{ m}{3}}\quad \mbox{ for }i=q+1,\dots, k.$$
This is possible, for $\lambda \binom{ m}{3}\leq (k-q)\rounddown{ m/3}$. We have
$$\dg_{G(i)}(\alpha)=3\mult_i(\alpha^3)\leq  m\quad \mbox{ for }i=q+1,\dots, k.$$
By the detachment lemma, there exists a hypergraph $G'$,  obtained by replacing the vertex $\alpha$ of $G$ by $ m$ new vertices  $\alpha_1,\dots,\alpha_{ m}$ in $G'$, replacing each $\alpha^3$-edge in $G$ by an $\alpha_i \alpha_j \alpha_\ell$-edge in $G'$ such that the following conditions hold (recall that $a\approx b$ means $a\in \{\rounddown{b}, \roundup{b}\}$).
\begin{align*}
    \dg_{G'(j)}(\alpha_i)&\approx\dfrac{\dg_{G(j)}(\alpha)}{ m}\leq 1 \quad \mbox{ for } i=1,\dots,  m, j=q+1,\dots, k;\\
    \mult_{G'}(\alpha_i \alpha_j \alpha_\ell)&=\dfrac{\mult_{G}(\alpha ^3)}{\binom{ m}{3}}=\lambda \quad \mbox{ for distinct } i,j,\ell\in \{1,\dots,  m\}.
\end{align*}
Thus, we obtain a $(k-q)$-coloring of $G'$. Since there is a one-to-one correspondence between $G'$ and  $\lambda  K_m^3$, we can color $\lambda  K_m^3\backslash F$ using the  $(k-q)$-coloring of $G'$. This together with the given $q$-coloring of $F$ lead to a $k$-coloring of $\lambda  K_m^3$. Applying Corollary \ref{n3mcor} completes the proof.
\qed

\section{Concluding Remarks and Open Problems} \label{remarksec}
Beside improving the bounds in Theorems \ref{weakerthm} and  \ref{evanstype}, and motivated by the Marcotte-Seymour theorem \cite{MR1073098}, one  may consider our  main problem without assuming that $ n\equiv 0\Mod 3$. 
\begin{problem} \label{mainprobgen}
Find conditions that ensure a coloring of $\lambda  K_m^3$ can be embedded in a coloring of $\lambda  K_n^3$ using the fewest possible number of colors.
\end{problem}
Similar problems are difficult even for the case of graphs \cite{MR3280683}. 
 Let $\chi'(\lambda  K_n^3)$ be the smallest number of colors needed to color $\lambda K_n^3$. 
\begin{theorem}
\begin{align*}
\chi'(\lambda  K_n^3)=
    \begin{cases}
    \lambda\dbinom{ n-1}{2} & \mbox {if }  n\equiv 0 \Mod 3,\\ 
    \\
    \lambda\dbinom{ n}{2} & \mbox {if }  n\equiv 2 \Mod 3,\\
    \\
    \dfrac{\lambda n( n-2)}{2} & \mbox {if }  n\equiv 4 \Mod 6,  \mbox {or if } n\equiv 1 \Mod 3, \lambda\equiv 0 \Mod 2, \\
    \\
    \dfrac{\lambda n^2-2\lambda  n+1}{2} & \mbox {if }  n\equiv 1 \Mod 6, \lambda\equiv 1 \Mod 2.
    \end{cases}
\end{align*}
\end{theorem}
\begin{proof}
We show that \begin{align} \label{chibound3}
    \chi'(\lambda  K_n^3)=\roundup{\frac{\lambda\dbinom{ n}{3}}{{\rounddown{\dfrac{ n}{3}}}}}.
\end{align}
The case $\lambda=1$ of \eqref{chibound3} was previously settled  in \cite{MR0416986}. Here, for convenience we provide a proof using hypergraph detachments. Our argument also works if we replace 3 by any $h< n$. It is clear that in any coloring of $\lambda  K_n^3$, each color class has at most $\rounddown{ n/3}$ edges. Therefore, 
    $\chi'(\lambda  K_n^3)\geq\roundup{\lambda\binom{ n}{3}/{\rounddown{\frac{ n}{3}}}}=:k$. 
To complete the proof, we find a $k$-coloring of $\lambda  K_n^3$. Let $G$ be a hypergraph whose edge set is $\lambda \binom{ n}{3}$ copies of $\{\alpha^3\}$, and let us color $G$ with colors $\{1,\dots, k\}$ such that 
$$\mult_{G(i)}(\alpha^3)\leq\rounddown{\frac{ n}{3}}\quad \mbox{ for }i=1,\dots, k.$$
This is possible, for $\lambda \binom{ n}{3}\leq k\rounddown{ n/3}$. By the detachment lemma, there exists a hypergraph $G'$,  obtained by replacing the vertex $\alpha$ of $G$ by  $\alpha_1,\dots,\alpha_{ n}$ in $G'$, replacing each $\alpha^3$-edge in $G$ by an $\alpha_i \alpha_j \alpha_\ell$-edge in $G'$ such that the following conditions hold.
\begin{align*}
    \dg_{G'(j)}(\alpha_i)&\approx\dfrac{\dg_{G(j)}(\alpha)}{ n}=\dfrac{3\mult_{G(j)}(\alpha^3)}{ n}\leq 1 \quad \mbox{ for } i=1,\dots,  n, j=1,\dots, k;\\
    \mult_{G'}(\alpha_i \alpha_j \alpha_\ell)&=\dfrac{\mult_{G}(\alpha ^3)}{\binom{ n}{3}}=\lambda \quad \mbox{ for distinct } i,j,\ell\in \{1,\dots,  n\}.
\end{align*}
This leads to a $k$-coloring of $G'\cong \lambda  K_n^3$. 

If $ n\equiv 0\Mod 3$, then $ n/3\in \mathbb Z$ and $k=3\lambda\binom{ n}{3}/ n=\lambda\binom{ n-1}{2}$.  If  $ n\equiv 2\Mod 3$, then $\rounddown{ n/3}=( n-2)/3$, and $k=3\lambda\binom{ n}{3}/( m-2)=\lambda\binom{ n}{2}$. Finally, if  $ n\equiv 1\Mod 3$, then $\rounddown{ n/3}=( n-1)/3$, and we have
\begin{align*}
k=\roundup{3\lambda\binom{ n}{3}/( n-1)}=\roundup{\frac{\lambda n( n-2)}{2}}=
\begin{cases}
\frac{\lambda n( n-2)}{2} & \mbox{ if }\lambda  n\equiv 0\Mod 2,\\
\frac{\lambda n( n-2)+1}{2} & \mbox{ if } \lambda \equiv  n\equiv 1\Mod 2.
\end{cases}    
\end{align*}

\end{proof}
An immediate consequence of the previous theorem is the following. 
\begin{corollary} \label{corisocoloring}
In any $\chi'(\lambda  K_n^3)$-coloring of $\lambda  K_n^3$, all color classes are isomorphic if and only if $ n\nequiv 1\Mod 6$ or $\lambda\equiv 0\Mod 2$. 
\end{corollary}
Corollary \ref{corisocoloring} reveals the very challenging case  when $ n\equiv 1 \Mod 6, \lambda\equiv 1 \Mod 2$, for in this case not all color classes are isomorphic.  Even in the other cases where all color classes are isomorphic but $ n\nequiv 0\Mod 3$, we encounter another issue. In order to extend  Theorem \ref{embedvialistthm}, to say, the case $ n\equiv 1 \Mod 3, \lambda\equiv 0 \Mod 2$, one has to decide which color classes have an isolated vertex in $\lambda K_m^3$.

Recall that for a function $\gamma$ from the edges of $G$ to $\mathcal P([k])$, a $\gamma$-coloring of $G$ is a coloring of $G$ in which the color of each edge $e$ of  $G$ is chosen from  $\gamma(e)$. Let  $\chi'_{\ell}(G)$ denote the smallest number $t$ such that $G$ has a $\gamma$-coloring whenever $|\gamma(e)|\geq t$ for each edge $e$ of $G$ (Here,  $k$ is sufficiently larger than $\chi'_{\ell}(G)$).  Let $f,g:[k]\rightarrow \mathbb N\cup \{0\}$ with $g(i)\leq \min \{f(i), |\{e\in G \ | \ i\in \gamma(e)|\}$ for $i\in [k]$. The following list coloring problem with restrictions  deserves investigating specially in connection with Theorem \ref{embedvialistthm} for the case where  $G=K_m$. 
\begin{problem} \label{listcolringquan}
Let $\gamma$ be a function  from the edges of $G$ to $\mathcal P([k])$ with $|\gamma(e)|\geq \chi'_{\ell}(G)$. Find conditions under which $G$ has a $\gamma$-coloring such that the following condition holds. 
\begin{align*}
    g(i)\leq |G(i)|\leq f(i) \mbox{ for } i\in [k].
\end{align*}
\end{problem}

Finally, we remark that there exists a symmetric layer-rainbow latin cubes of order $n$  if and only if $n\equiv 0,2 \Mod 3$ (with two exceptions;  $n$ can be 1, but $n$ cannot be 3) \cite{doi:10.1137/22M1494488}. Hence, it seems reasonable to conjecture the following.
\begin{conjecture}
    There exists a constant $c\geq 3$ such that the following holds.  If $n\geq cm$, then an $m\times m\times m$  symmetric layer-rainbow latin box on $n$ symbols can be embedded in a symmetric layer-rainbow latin cube of order $n$ if and only if $n\equiv 0,2 \Mod 3$. 
\end{conjecture}

\section{Acknowledgment} We wish to thank the anonymous referees for their constructive criticism.

\bibliographystyle{plain}

\end{document}